\documentclass[12pt]{article} 

\usepackage[reqno]{amsmath}
\usepackage{amssymb,amsthm}

\newtheorem{theorem}{Theorem}
\newtheorem{lemma}{Lemma}

\newcommand{\ph}{\varphi}
\newcommand{\hb}{\frac{1}{2}}
\newcommand{\h}{1/2}
\newcommand{\eps}{\varepsilon}
\newcommand{\del}{\delta}
\newcommand{\zeile}{\vspace{\baselineskip}}

\newcommand{\N}{\mathbb{N}}

\newcommand{\Z}{\mathbb{Z}}
\newcommand{\R}{\mathbb{R}}
\newcommand{\C}{\mathbb{C}}

\allowdisplaybreaks

\begin{document}

\title{A new bound for the large sieve inequality\\ with power
    moduli}

\author{Karin Halupczok}
\maketitle


\begin{abstract}
We give a new bound for the large sieve inequality with power moduli
$q^{k}$ that is uniform in $k$. The proof uses a new theorem due
to T.\ Wooley from his work on efficient congruencing.
\end{abstract}

Kewords: Large sieve inequality, Powers, Weyl sums

MSC 2010: 11N35, 11L15

\section{Introduction}

Let $\{v_{n}\}$ denote a sequence of complex numbers, $M,N,k\in\N$,
and let $Q$ be a real number $\geq 1$. 
We write $e(\alpha):= \exp(2\pi i \alpha)$ for $\alpha\in\R$.

The large sieve inequality with power
moduli aims to give upper bounds for the sum
\[
   \Sigma_{Q,N,k}:=   \sum_{q\leq Q} \sum_{\substack{1\leq a\leq
       q^{k}\\ \gcd(a,q)=1}}  \Big| \sum_{M<n\leq M+N} 
       v_n    e\Big(\frac{a}{q^{k}}n\Big) \Big|^2.
\]

It is known that an application of the
standard large sieve inequality gives the upper bounds
\begin{equation}
\label{e:ls}
   \Sigma_{Q,N,k}\ll_{k} (N+Q^{2k}) |v|^{2} \;\text{ and }\;
   \Sigma_{Q,N,k}\ll_{k} (QN+Q^{k+1}) |v|^{2},
\end{equation}
where $|v|^{2}:=  \sum_{M<n\leq M+N} |v_n|^2$,
and it is conjectured by L.\ Zhao in
\cite{c1} that the upper bound
\begin{equation}
\label{e:conj}
  \Sigma_{Q,N,k} \ll_{k,\eps} |v|^{2} (N+Q^{k+1})(NQ)^{\eps}
\end{equation}
should hold. 

The bounds (\ref{e:ls}) verify the conjecture 
for $Q\leq N^{1/(2k)}$ and $Q\geq N^{1/k}$, so the problem is to prove
it in the range 
\begin{equation*}
    N^{1/(2k)}\leq Q\leq N^{1/k} \Leftrightarrow Q^{k}\leq N \leq Q^{2k}.
\end{equation*}

Especially the cases for small $k$, namely $k=2,3$
are of interest and were considered in the
papers \cite{c2},\cite{c2a} and \cite{c1}.
In this paper we investigate 
the problem uniform in $k$. The following
nontrivial bounds are known in this case.

L.\ Zhao showed in \cite{c1} the bound
\begin{equation}
\label{eq:Z}
    \Sigma_{Q,N,k} \ll_{k,\eps} |v|^{2} (Q^{k+1} + 
     (NQ^{1-1/\kappa} + N^{1-1/\kappa}Q^{1+k/\kappa}) N^{\eps}),
\end{equation}
where $\kappa:=2^{k-1}$.

In \cite{c2a}, it was shown by S.\ Baier and L.\ Zhao that
\begin{equation}
\label{eq:BZ}    \Sigma_{Q,N,k} \ll_{k,\eps} |v|^{2} (Q^{k+1} + N +
    N^{1/2+\eps}Q^{k})(\log\log 10 NQ)^{k+1}
\end{equation}
holds, which improves Zhao's bound (\ref{eq:Z}) for
$Q \ll N^{(\kappa-2)/(2(k-1)\kappa -2k)-\eps}$.

In this paper we prove the following result:

\vspace{1ex}
\begin{theorem}
\label{t}
 Let $\del:=(2k(k-1))^{-1}$.
 Then we have the bound
  \begin{equation*}
     \Sigma_{Q,N,k} \ll_{k,\eps}
        |v|^{2} (NQ)^{\eps}(Q^{k+1} + Q^{1-\del} N + Q^{1+k\del} N^{1-\del}).
  \end{equation*}
\end{theorem}

This bound improves the bound (\ref{eq:Z}) for all $k$
sufficiently large, and the bound (\ref{eq:BZ}) for
$ Q^{k}\leq  N \leq Q^{2k-2+2\del}$ and all $k\geq 3$, but it does not
confirm any case of Zhao's conjecture (\ref{e:conj}), too.
Further, the result is not sufficient to give 
an improvement of the bound in \cite{c2a} for $k=3$, but comes
near to it. 

\zeile
\textbf{Notation.} In the following,
we suppress the dependence of the implicit constants
on $k$ or $\eps$ in our estimates and write simply
$\ll$ for $\ll_{k,\eps}$. The small value $\eps>0$ may
depend on $k$ and may change its value during calculation.
The symbol $\|\alpha\|$ means the distance of
$\alpha$ to the nearest integer, and by $\{\alpha\}:=\alpha-[\alpha]$
we denote the fractional part of $\alpha$, and by $[\alpha]$
the largest integer smaller or equal to $\alpha$.

\section{Lemmas}

We make use of the following version of the large sieve inequality. 

\begin{lemma}
\label{l1}
  Let $S$ denote a finite set of positive integers, 
  $M,N\in\Z$ and let $\{v_n\}$ be a complex sequence. Further let
  \[
    \mathcal{F}:=\{(a,q)\in\Z^{2};\; q\in S,\; 0<a<q,\; \gcd(a,q)=1\}.
  \]
  Then
  \begin{multline}
   \label{eq:1}
      \sum_{(a,q)\in\mathcal{F}}
         \Big| \sum_{M<n\leq M+N} 
       v_n e\Big(\frac{a}{q}n\Big) \Big|^2  \\
     \leq  \sum_{M<n\leq M+N} |v_n|^2 \Big( 4\sum_{q\in S} q + 
       \max_{(b,r)\in\mathcal{F}} \int_{1/N}^{\h}
          \#\mathcal{F}_{b,r}(x)  \frac{dx}{x^2}  \Big),
  \end{multline}
  \end{lemma}
  where 
  \[
    \mathcal{F}_{b,r}(x):=\Big\{ (a,q)\in\mathcal{F};\; 
     \Big|\frac{a}{q} - \frac{b}{r} \Big|\leq x \Big\}.
  \]

\textbf{Proof:}
We use Halasz-Montgomery's inequality
\[
    \sum_{r\leq R} |\langle v,\ph_r\rangle|^2 \leq
      |v|^2 \cdot \max_{r\leq R} 
     \sum_{s\leq R} |\langle \ph_r,\ph_s \rangle|
\]
that holds for any sequence $\{\ph_r\}$ of vectors of $\C^{N}$, 
and where $|v|^2=\langle v,v \rangle$, and 
$\langle \cdot , \cdot\rangle$ is the standard scalar
product on $\C^N$.

So the left hand side of (\ref{eq:1}) is
\begin{align*}
     \sum_{(a,q)\in\mathcal{F}} &\Big| \sum_{M<n\leq M+N} 
       v_n e\Big(\frac{a}{q}n\Big) \Big|^2  \\
     &\leq |v|^2 \max_{(b,r)\in \mathcal{F}} 
       \sum_{(a,q)\in \mathcal{F}}
    \Big| \sum_{M<n\leq M+N} e\Big(\frac{a}{q}n\Big)
     e\Big(-\frac{b}{r}n\Big) \Big|   \\
   &\leq  |v|^2  \max_{(b,r)\in \mathcal{F}} \sum_{(a,q)\in \mathcal{F}}
      \min\Big(N,\Big\|\frac{a}{q}-\frac{b}{r}\Big\|^{-1}\Big).
\end{align*}

Now we have to estimate

\begin{equation}
   \label{eq:2}
   \max_{(b,r)\in \mathcal{F}} \sum_{(a,q)\in \mathcal{F}}
      \min\Big(N,\Big\|\frac{a}{q}-\frac{b}{r}\Big\|^{-1}\Big).
\end{equation}

For this, fix $(b,r)\in \mathcal{F}$. For $\Delta>0$ write
\[
   P(\Delta):= \# \mathcal{F}_{b,r}(\Delta).
\]

Let $\Delta_0:=\frac{1}{N}$ and for $L\in\N$ let 
$h:=(\hb-\frac{1}{N})L^{-1}$. 
Now let $\Delta_i:=\frac{1}{N}+hi$, 
so $\Delta_L=\frac{1}{2}$.
Since $\|\alpha\|=\min\{|\alpha|,1-|\alpha|\}$ for $-1<\alpha<1$,
we have
\begin{align*} 
   &\sum_{(a,q)\in \mathcal{F}} 
   \min\Big(N,\Big\|\frac{a}{q}-\frac{b}{r}\Big\|^{-1}\Big) \\
   &\leq 2NP\Big(\frac{1}{N}\Big) + 2\sum_{0\leq i < L} 
   \sum_{\substack{(a,q)\in \mathcal{F}\\ \Delta_i < |a/q-b/r| 
    \leq \Delta_{i+1}}} \frac{1}{\Delta_i} \\
   &=2NP\Big(\frac{1}{N}\Big) + 2\sum_{0\leq i < L} \frac{1}{\Delta_i} 
    (P(\Delta_{i+1})-P(\Delta_i)) \\
   &=2\sum_{0\leq i < L}
    \Big(\frac{1}{\Delta_i}-\frac{1}{\Delta_{i+1}}\Big)
   P(\Delta_{i+1}) + \frac{2}{\Delta_{L}} P(\Delta_L).
\end{align*}

The last summand is $\leq 4\sum_{q\in S} q$, and
the sum over $i$ approximates the Riemann-Stieltjes-integral
\begin{equation}
\label{RSI}
   \int_{1/N}^{\h} P(x) dg(x) \text{ with } g(x)=-\frac{1}{x},
\end{equation}
if $L\to\infty$. Therefore the sum over $(a,q)\in \mathcal{F}$ in (\ref{eq:2})
is at most as large as the integral (\ref{RSI}), plus $4\sum_{q\in S} q$.

Since $g$ is 
continuously differentiable
on $[\frac{1}{N},\frac{1}{2}]$ 
and since $P$ is  Riemann-integrable, the integral (\ref{RSI}) equals
\[
   \int_{1/N}^{\h} P(x) g'(x) dx
   = \int_{1/N}^{\h} P(x) \frac{dx}{x^2} 
= \int_{1/N}^{\h} \# \mathcal{F}_{b,r}(x)
     \frac{dx}{x^2}.
\]
This was to be shown. \qed

\zeile\zeile
Further we use the following estimate for the 
exponential sum occurring in the proof of Theorem \ref{t}.

\begin{lemma}
  \label{l2}
  Let $f(x):=\alpha x^{k}\in\R[x]$ be a monomial of degree $k\geq 2$,
  and $S_{Q}:=\sum_{Q<q\leq 2Q} e(f(q))$, $\del:=(2k(k-1))^{-1}$. 
  Then 
  \[
      S_{Q} \ll Q^{1+\eps} \Big(Q^{-1} + Q^{-k}\sum_{1\leq v\leq Q} 
       \min(Q^{k}v^{-1},\|v\alpha\|^{-1})\Big)^{\del}.
  \]
\end{lemma}

\textbf{Proof:}

Suppose that $a,q\in\Z$ with $(a,q)=1$ and $|q\alpha-a|\leq q^{-1}$.

We apply Theorem 1.5 in T.\ Wooley's article \cite{c1} on
efficient congruencing and obtain
\[
   S_{Q} \ll Q^{1+\eps} (q^{-1}+Q^{-1}+qQ^{-k})^{\del}.
\]
By a standard transference principle (see Ex.\ 2 of section 2.8
in Vaughan's book \cite{c7}), this implies that
\begin{equation}
\label{eq:sq}
   S_{Q} \ll Q^{1+\eps} \Big( (v+Q^{k}|v\alpha-u|)^{-1} + Q^{-1} +
   (v+Q^{k}|v\alpha-u|)Q^{-k} \Big)^{\del}
\end{equation}
for any integers $u,v\in\Z$ with $(u,v)=1$ and $|v\alpha-u|\leq v^{-1}$.

Now by Dirichlet's Approximation Theorem, there exist such integers
$u,v$ with $1\leq v\leq Q^{k-1}$ and $|v\alpha-u|\leq Q^{1-k}$, for these
\[
   (v+Q^{k}|v\alpha-u|)Q^{-k} \ll (Q^{k-1}+Q)Q^{-k} \ll Q^{-1}
\]
holds. Further we get
\[
   (v+Q^{k}|v\alpha-u|)^{-1} \ll Q^{-k}
   \min(Q^{k}v^{-1},|v\alpha-u|^{-1}).
\]

Now if $v>Q$, this expression is again $\ll Q^{-1}$.
If otherwise $1\leq v\leq Q$, it is bounded by
\[
   Q^{-k} \sum_{1\leq v\leq Q}\min(Q^{k}v^{-1},\|v\alpha\|^{-1}),
\]
since $|v\alpha-u|\geq\|v\alpha\|$.

Hence, these estimates included in (\ref{eq:sq}) show 
the assertion. 

\qed

\zeile
\begin{lemma}
  \label{l3}
  Let $X,Y,\alpha\in\R$, $X,Y\geq 1$, and
  $a,q\in\Z$, $\gcd(a,q)=1$, with $|q\alpha-a|\leq q^{-1}$.
  Then
  \[
     \sum_{v\leq X} \min\Big(XYv^{-1}, \|\alpha v\|^{-1} \Big)
     \ll XY (q^{-1}+Y^{-1}+q(XY)^{-1})\log (2Xq).
  \]
\end{lemma}

This is Lemma 2.2 of \cite{c7}. \qed

\section{Proof of Theorem \ref{t}}

Let $k\in \N$ with $k\geq 2$, let 
$Q\geq 1$ and assume that the integer $N$ is in the range
$Q^{k}\leq N\leq Q^{2k}$.

We apply Lemma \ref{l1} with
\[
   \mathcal{F}:=\{(a,q^{k})\in \Z^{2};\; 
   Q< q\leq 2Q,\; 0<a<q^{k},\; \gcd(a,q)=1  \},
\]
which shows that
\[
   \Sigma_{Q,N,k} \ll |v|^{2}Q^{\eps} \Big(Q^{k+1} + 
      \max_{(b,r^{k})\in\mathcal{F}} \int_{1/N}^{1/2}
      \#\mathcal{F}_{b,r^{k}}(x) \frac{dx}{x^{2}} \Big),
\]
since we have the admissible error $\sum_{q\leq Q}q^{k}\ll Q^{k+1}$.

Now we aim to give an estimate for
\begin{equation*}
    \max_{(b,r^{k})\in\mathcal{F}} \int_{1/N}^{1/2}
      \#\mathcal{F}_{b,r^{k}}(x) \frac{dx}{x^{2}}
\end{equation*}

The integrand counts for fixed $(b,r^{k})\in\mathcal{F}$ 
all $(a,q^{k})\in \mathcal{F}$ with
\[
   \Big| \frac{a}{q^{k}} - \frac{b}{r^{k}} \Big| \leq x.
\]

So for fixed $Q<q\leq 2Q$, we count every $a$ with
\[
   \frac{|ar^{k}-bq^{k}|}{2r^{k}Q^{k}x}\leq \hb.
\]

Now we use the Fourier analytic method from the papers 
\cite{c2},\cite{c2a} and \cite{c1} by Baier
and Zhao. For this, consider the function 
\[
   \phi(x):= \Big(\frac{\sin \pi x}{2x}\Big)^{2}, \qquad\qquad
\phi(0):= \lim_{x\to 0}\phi(x)  =\frac{\pi^{2}}{4}.
\]
Then $\phi(x)\geq 1$ for $|x|\leq 1/2$, and the Fourier transform of
$\phi$ is 
\begin{equation*}
    \hat{\phi}(s) = \frac{\pi^{2}}{4} \max\{1-|s|,0 \}.
\end{equation*}

For fixed $q$, we get for the number of 
corresponding $a$ the estimate
\[
     \sum_{a,(a,q)\in\mathcal{A}_{b,r}(x)} 1
     \leq  \sum_{a\in \Z}
     \phi\Big(\frac{ar^{k}-bq^{k}}{2r^{k}Q^{k}x}\Big)
     =  \sum_{a\in\Z} \int_{-\infty}^{\infty}  
      \phi\Big(\frac{sr^{k}-bq^{k}}{2r^{k}Q^{k}x}\Big) e(as)ds, 
\]
where we applied in the last step Poisson's summation formula.
Summing up over $q$ and a linear transformation gives
\begin{multline*}
   \sum_{Q<q\leq 2Q}\sum_{a\in\Z} \int_{-\infty}^{\infty}
          \phi(v) e\Big(ab \frac{q^{k}}{r^{k}} \Big)
        e(2Q^{k}xav) 2Q^{k}x dv \\
   = \sum_{|a|\leq B} \hat{\phi}\Big(\frac{a}{B}\Big)B^{-1} 
     \sum_{Q<q\leq 2Q} e\Big(ab\frac{q^{k}}{r^{k}}\Big),
\end{multline*}
where we have set $B:=(2Q^{k}x)^{-1}$, and we may assume w.l.o.g.\
that $B\geq 1$.

We separate the summand with $a=0$ and get
\begin{equation*}
   \ll Q^{k+1}x + B^{-1} \sum_{1\leq a\leq B}
    \Big| \sum_{Q<q\leq 2Q} e \Big( \frac{abq^{k}}{r^{k}}\Big) \Big|.
\end{equation*}

The separated term $Q^{k+1}x$ leads again to the admissible 
contribution
\[
    \int_{1/N}^{\h} Q^{k+1} \frac{dx}{x} \ll Q^{k+1+\eps}.
\]

Consider the monomial $f(q):=\frac{ab}{r^{k}}q^{k}$ of 
degree $k$ in $q$ and coefficient $\alpha:=\frac{ab}{r^{k}}\neq 0$.
It remains to give a good upper bound for the expression
\begin{equation}
\label{eq:e6}
   \int_{1/N}^{\h} B^{-1} \sum_{1\leq a\leq B}   
     \Big| \sum_{Q\leq q<2Q} e (f(q)) \Big| \frac{dx}{x^{2}}.
\end{equation}
Denote by $S_{Q}$ the occurring exponential sum
\[
   S_{Q}:= \sum_{Q<q\leq 2Q} e(f(q)).
\]
By Lemma \ref{l2}, we have 
\[
    S_{Q} \ll Q^{1+\eps} \Big(Q^{-1} + Q^{-k}\sum_{1\leq v\leq Q} 
     \min(Q^{k}v^{-1},\|v\alpha\|^{-1})\Big)^{\del}.
\]
The summand $Q^{-1}$ in big parantheses provides already the contribution
\begin{equation}
\label{eq:P2}
    \int_{1/N}^{\h}  Q^{1-\del+\eps} \frac{dx}{x^{2}} \ll Q^{1-\del+\eps}N
\end{equation}
to (\ref{eq:e6}),
and it remains to consider the term with the sum over $v$.

We estimate its contribution to $S_{Q}$ as follows using 
H\"{o}lder's inequality and Lemma \ref{l3}. We have
\begin{align*}
  &Q^{1+\eps-k\del} \sum_{a\leq B} \Big(\sum_{v\leq Q}
  \min\Big(Q^{k}v^{-1},\Big\|\frac{ab}{r^{k}} v
  \Big\|^{-1}\Big)\Big)^{\del} \\
 \ll &Q^{1+\eps-k\del} B^{1-\del}\Big( \sum_{\ell\leq BQ} d(\ell) 
   \min\Big(BQ^{k}\ell^{-1},\Big\|\frac{b}{r^{k}}\ell \Big\|^{-1}\Big) 
    \Big)^{\del} \\
 \ll &Q^{1+\eps-k\del} B^{1-\del} \Big((BQ^{k})^{1+\eps}(r^{-k}
    +Q^{1-k}+r^{k}(BQ^{k})^{-1}) \Big)^{\del}  \\
 \ll &BQ^{1+\eps}(Q^{1-k}+ B^{-1})^{\del}.
\end{align*}

The contribution to (\ref{eq:e6}) becomes
\begin{align*}
   &\ll  Q^{1+\eps+(1-k)\del} N + Q^{1+\eps}\int_{1/N}^{\h}
   B^{-\del}\frac{dx}{x^{2}} \\ &\ll Q^{1-(k-1)\del+\eps} N 
   + Q^{1+\eps} \int_{1/N}^{\h} Q^{k\del}x^{\del} \frac{dx}{x^{2}}
   \\ &\ll Q^{1-(k-1)\del+\eps} N + Q^{1+k\del+\eps} N^{1-\del}.
\end{align*}
The first term can be estimated by the bound
(\ref{eq:P2}), since $k\geq 2$. 
We obtain the stated bound of Theorem \ref{t}. 
\qed


\section*{Acknowledgments}
Thanks to T.\ Wooley for the suggestion which helped
to improve Lemma \ref{l2} and a former result.
Thanks also to S.\ Baier for comments concerning this paper.


\end{document}